\documentclass[12pt]{article} 
\usepackage{latexsym}
\usepackage{amssymb}
\usepackage{euscript}

\let\cal\mathcal
\usepackage{latexsym}
\usepackage{amssymb}
\usepackage{euscript}
\usepackage[dvips]{graphics}
\usepackage{epsf}

\let\cal=\mathcal      

\def\mcc{M\raise.5ex\hbox{c}C}
\def\mccarthy{M\raise.5ex\hbox{c}Carthy}

\def\eg{{\it e.g. }}
\def\ie{{\it i.e. }}

\def\h{{\cal H}}

\def\K{{\cal K}}

\def\N{{\cal N}}
 
\def\LL{{\cal L}}

\def\m{Mult}

\def\a{\alpha}
\def\l{\lambda}

\def\vare{\varepsilon}

\let\i=\infty

\def\la{\langle}
\def\ra{\rangle}
\def\={\ = \ }

\def\C{\mathbb C}
\def\R{\mathbb R}

\def\D{\mathbb D}

\def\be{\setcounter{equation}{\value{theorem}} \begin{equation}}
\def\ee{\end{equation} \addtocounter{theorem}{1}}
\def\beq{\begin{eqnarray*}}
\def\eeq{\end{eqnarray*}}
\def\se{\setcounter{equation}{\value{theorem}}} 
\def\att{\addtocounter{theorem}{1}}
\def\vs{\vskip 5pt}
\def\bs{\vskip 12pt}

\def\bp{{\sc Proof: }}
\def\ep{{}{\hfill $\Box$} \vskip 5pt \par}

\def\obc{{\sc Proof of Claim: }}

\def\bl{\begin{lemma}}
\def\el{\end{lemma}}
\def\bt{\begin{theorem}}
\def\et{\end{theorem}}
\def\bprop{\begin{prop}}
\def\eprop{\end{prop}}
\def\bd{\begin{definition}}
\def\ed{\end{definition}}
\def\br{\begin{remark}}
\def\er{\end{remark}}
\def\bexer{\begin{exercise}}
\def\eexer{\end{exercise}}

\newtheorem{theorem}{Theorem}[section]
\newtheorem{prop}[theorem]{Proposition}
\newtheorem{lemma}[theorem]{Lemma}
\newtheorem{cor}[theorem]{Corollary}

\newtheorem{conjecture}[theorem]{Conjecture}

\def\w{{\bf w}}
\def\wa{{\bf w}^\alpha}
\def\hw{{{\cal H}_\w}}
\def\hwf{{{\cal H}_{\w^\flat}}}
\def\ha{{{\cal H}_\alpha}}
\def\a{\alpha}
\def\s{\sigma}
\def\O{\Omega}
\def\Oh{\O_{1/2}}
\def\Oo{{\O_0}}
\def\D{\cal D}
\def\N{\mathbb N}
\def\f{\flat}
\def\phip{{\phi^\prime}}

\begin{document}
\setlength{\baselineskip}{21pt}
\title{Hilbert spaces of Dirichlet Series}
\author{
John E. M\raise.5ex\hbox{c}Carthy
\thanks{Partially supported by National Science Foundation Grant
DMS 0070639}\\
Washington University\\
St. Louis, Missouri 63130}
\date{December 7 2001}

\bibliographystyle{plain}

\maketitle
\begin{abstract}
We consider various Hilbert spaces of Dirichlet series whose norms
are given by weighted $\ell^2$ norms of the Dirichlet coefficients.
We characterize the multiplier algebras for some of these spaces.
\end{abstract}

\baselineskip = 18pt

\setcounter{section}{-1}
\section{Introduction}
Let $\w \, = \, \{w_n\}_{n=n_0}^\i$ be a sequence of positive
numbers. In this paper we are concerned with Hilbert spaces of
functions representable by Dirichlet series:
\be
\hw \ =\ \left\{ f(s) \, = \, \sum_{n=n_0}^\i a_n n^{-s} \ \Big| \
\| f \|^2_{\hw} \, := \, \sum_{n=n_0}^\i |a_n|^2 w_n \, < \, \i \right\}
.
\ee

The prototypical case, where $w_n \equiv 1,\ n \geq 1$, was first 
studied in a
beautiful paper by H.~Hedenmalm, P.~Lindqvist and K.~Seip
\cite{hls}. Among other results, they characterized the multipliers
of the space (see Theorem~\ref{thmb2} below for a statement of their
result).

One purpose of this paper is to consider the scale of spaces
obtained from the weight sequences $\wa$, defined for $n \geq 2$ by
\be
\label{eq02}
w^\alpha_n \= ( \log n )^\alpha .
\ee
For brevity, we shall write $\ha$ for the space $\h_{\wa}$;
specifically 
$$
\ha \= \left\{ f(s) \, = \, \sum_{n=2}^\i a_n n^{-s} \ \Big| \
 \sum_{n=2}^\i |a_n|^2 (\log n)^\alpha \, < \, \i \right\}
.$$
(When $\alpha = 0$, it is more natural to let $n_0 =1$ and to include the constant functions
in $\h_0$. 
It is not essential to any of the issues we discuss here).
\bs

Before going further, let us remind the reader of some basic facts about Dirichlet series.
A nice treatment can be found in Titchmarsh's book \cite{tit32}.
We shall follow the convention of writing the complex variable $s = \s + it$.
A Dirichlet series is a series of the form 
\be
\label{eq3}
\sum_{n=1}^\i a_n n^{-s} .
\ee
Such a series may converge for no values of $s$; if it converges for
any particular $s_0$, then it converges for all $s$ with $\s > \Re(s_0)$. Therefore the largest open set in which a 
series (\ref{eq3}) converges is a half-plane (at what points on the boundary of the 
half-plane the series converges is, in general, a delicate question).
Let us adopt the notation, for $\rho$ a real number, $\O_\rho$ is
the half-plane
$$
\O_\rho \= \{ s \, \in \, \C \ | \ \sigma > \rho \}  .
$$
Let 
$$
\s_c = \inf \{ \Re(s) \, : \,  
\sum_{n=1}^\i a_n n^{-s} \ {\rm converges} \}; 
$$
this is called the 
\emph{abscissa of convergence} of the series.
The largest domain of convergence of the series is $\O_{\s_c}$.
There are three
other abscissae associated with the series (\ref{eq3}) which we shall need.
The first is the 
\emph{abscissa of absolute convergence}, $\s_a$, defined by
$$
\s_a = \inf \{ \Re(s) \, : \,  
\sum_{n=1}^\i a_n n^{-s} \ {\rm converges \ absolutely} \}.
$$
Obviously $\s_a \geq \s_c$; it is straightforward 
that $\s_a \leq \s_c + 1$, because if $\sum_{n=1}^\i a_n n^{-s}$
converges, then $|a_n| n^{-\s} = o(1)$.
The second
is the \emph{abscissa of boundedness} $\s_b$, defined by
$$
\s_b = \inf \{ \rho \, : \,  
\sum_{n=1}^\i a_n n^{-s} \ {\rm converges\ to\ a\ bounded\ function\
in\ }\O_\rho \}. 
$$
The third abscissa is 
the \emph{abscissa of uniform convergence} $\s_u$, defined by
$$
\s_u = \inf \{ \rho \, : \,
\sum_{n=1}^\i a_n n^{-s} \ {\rm converges\ uniformly\
in\ }\O_\rho \}.
$$
Clearly $\s_u \geq \s_b$;
H.~Bohr proved that $\s_u = \s_b$, and 
$\s_a \leq \s_b + \frac 12$
\cite{boh13}. Note that the series 
does not necessarily define
a bounded function in $\O_{\s_b}$, but the function it represents 
is bounded in all \emph{strictly smaller} half-planes.
If all the coefficients $a_n$ are positive, then all three of $\s_c,
\s_b, \s_a$ coincide.

We shall let $\D$ denote the set of functions that can be 
represented in some half-plane by
a Dirichlet series. 

Let $f(s)$ be holomorphic in the half-plane $\O_\rho$. Let $\vare >
0$. A real number $\tau$ is called an \emph{$\vare$ translation
number of $f$} if
$$
\sup_{s \in \O_\rho} | f(s+ i \tau) - f(s) | \leq \vare .
$$
The function $f(s)$ is called \emph{uniformly almost periodic}
in the half-plane $\O_\rho$ if, for every $\vare > 0$, there exists
a positive real number $M$   such that every interval in $\R$ of
length $M$ contains at least one $\vare$ translation
number of $f$. 

We shall need the following theorem. The proof 
can
be found in \cite[p. 144]{bes}.
\bt
\label{thma2}
Suppose that $f(s)$ is represented by a Dirichlet series that
converges uniformly in the half-plane $\O_\rho$. Then $f$ is
uniformly almost periodic in $\O_\rho$.
\et

\bs

Returning to the spaces $\ha$, it follows from the Cauchy-Schwarz inequality
that any function in any space $\ha$ has $\s_a \leq \frac 12$. Moreover, for all
$\vare >0$, the function\footnote{
We use $\zeta$ to denote the Riemann zeta-function: $\zeta(s) =
\sum_{n=1}^\i n^{-s}$.
}
$$
\zeta( \frac 12 + \vare + s) -1 \= \sum_{n=2}^\i \frac{1}{n^{\frac 12 + \vare}} n^{-s}
$$
is in every $\ha$ and has a pole at $\frac 12 - \vare$, so the largest common
domain of analyticity of the functions in any $\ha$
is $\Oh$. The reproducing kernel for $\hw$ is
\be
\label{eqa35}
k(s,u) \= \sum_{n=n_0}^\i \frac{1}{w_n} n^{-s - \bar u};
\ee
for the spaces $\ha$ this is essentially a fractional derivative or integral of the $\zeta$
function at $s + \bar u$.

\bs
We claim that the scale of spaces $\ha$ is in many ways analogous
to the scale of spaces of holomorphic functions in the unit disk defined by
\be
\label{eq4}
\K_\alpha \= \left\{ g(z) \, = \, \sum_{n=0}^\i a_n z^n \ \Big| \
\sum_{n=0}^\i |a_n|^2 (n+1)^\alpha \, < \, \i \right\}.
\ee
Thus the space $\h_0$ corresponds to the Hardy space, the space
$\h_{-1}$ corresponds to the Bergman space, and the space $\h_1$
corresponds to the Dirichlet space in the setting of the disk.

Notice first that differentiation is a unitary map from $\ha$ to $\h_{\alpha -2}$, 
and fractional differentiation/integration maps any $\ha$ to any $\h_\beta$.
Although this is not exactly true for the spaces $\K_\alpha$ as we defined them, it is essentially
true. It fails only because (a) differentiation annihilates the constant functions, and 
(b) the choice of norms in (\ref{eq4}) does not quite render differentiation
isometric. For example, the norm of $z^n$ in the Dirichlet space
is 
$\sqrt{n+1}$, and the norm of $nz^{n-1}$ in the Bergman space 
is $\sqrt{n}$. Asymptotically these agree, and it is true as function spaces that the Dirichlet space 
is the set of functions whose derivatives are in the Bergman space.

The main reason for our analogy is the Plancherel formula (Theorem~\ref{thmb1} below)
which gives the norm of a function in $\ha$, for $\alpha > 0$, in terms of a weighted area
integral. As $\alpha$ tends to $0$, the measure tends to a line integral on the boundary.
Finally, the multiplier algebras of the spaces $\ha$ are the same for
all values of $\alpha \leq 0$, and shrink as $\alpha$ becomes positive
(Theorems~\ref{thmb2} and \ref{thmc1}). 

Why study these spaces? We offer two considerations. First, 
a Plancherel formula calculating the norm of a function in
two different ways
frequently has serendipitous consequences (see \eg \cite{hls}). Second, the kernel
functions for these spaces are formed from the zeta function and
its relatives (\ref{eqa35}), and one
may hope that studying the spaces will shed further light on the
functions.

\bs
In Section~\ref{secpick}, we consider a different choice of weight
sequence that gives rise to a space of Dirichlet series with
reproducing kernel
$$
k(s,u) \= \frac{1}{2 - \zeta(s + \bar u)} .
$$
This has the interesting property of being a complete Pick kernel,
and several properties flow from this (see \eg \cite{ampi} for a
discussion). We show that, unlike for the spaces $\h_\alpha$, the
multipliers do not extend to be analytic on a larger domain 
than the domain of analyticity of the Hilbert space. 

\section{Bergman-like Spaces}
\label{secber}

Throughout this section, $\mu$ will be a positive Radon measure on
$[0,\i)$ for which 
\be
\label{eq11}
\int_0^\i n_0^{-2 \sigma} d\mu(\sigma) \= \int_0^\i e^{-2 (\log n_0) \,
\sigma} d\mu(\sigma) \ < \ \i
\ee
for some positive integer $n_0$. We also assume that 
\be
\label{eq12}
0\ {\rm is\ in\ the\
support\ of\ }\mu.
\ee
We define $w_n$, for $n \geq n_0$, by
\be
\label{eqb3}
w_n \ := \ 
\int_0^\i n^{-2 \sigma} d\mu(\sigma) .
\ee

Letting $\mu_\alpha$ be the measure
$$
d\mu_\alpha (\s) \= \frac{2^{-\alpha}}{\Gamma(-\a)} \, \s^{-1-\a} d\s
$$
gives the weights from (\ref{eq02}) for $\alpha < 0$: 
$$
\int_0^\i n^{-2 \sigma} d\mu_\alpha(\sigma) \= (\log n)^\a, \qquad n
\geq 2, \ \a < 0.
$$
We let $\mu_0$ be the unit point mass at $0$, which has all of its
moments equal to $1$.

For any measure satisfying (\ref{eq11}) and (\ref{eq12}), the moments
$w_n$ are a decreasing sequence that decays more slowly than any
negative power of $n$: for all $\vare > 0$, there exists $c >0 $
so that
\be
\label{eq13}
w_n \ > \  c\, n^{-\vare}.
\ee
Therefore, every space $\hw$ consists of functions analytic in $\Oh$,
and contains functions that are not analytically extendable to any
larger domain. Nonetheless there is a dense subspace of functions in $\hw$ 
whose norms can be obtained by evaluating suitable integrals over the
larger half-plane $\Oo$. For the case of $\h_0$, the following
theorem is due to F.~Carlson \cite{car22}.
Note that we do not assume that either side of (\ref{eqb4}) is finite.

\bt
\label{thmb1}
Let $f(s) = \sum_{n=n_0}^\i a_n n^{-s}$ be a function in $\D$ that has
$\s_b =0$. Let $\mu$ satisfy (\ref{eq11}) and (\ref{eq12}), and let $w_n$
be given by (\ref{eqb3}). Then 
\be
\label{eqb4}
\sum_{n=n_0}^\i |a_n|^2 w_n
\= \lim_{c \to 0^+} \, \lim_{T \to \i}\,  \frac{1}{2T}
\int_{-T}^T dt \int_0^{\i} d\mu(\sigma) | f(s + c)|^2 .
\ee
\et
\bp
Fix $ 0 < c < 1$, and let $\vare >0$. 
Let $\delta $ be given by
$$
\delta \= \frac{\vare}{(1 + \mu[c, \frac 1c + c])\, ( 1 + \| f
\|_{\O_c} )} .
$$
By Bohr's theorem, the Dirichlet series 
of $f$ converges 
uniformly in $\overline{\O_c}$, 
so there exists $N$ so that
$$
| \sum_{n=n_0}^{N^\prime} a_n n^{-s} - f(s) | \ < \ \delta \qquad
\forall \, s \, \in \, \overline{\O_c}, \ \forall \,  N^\prime \, \geq
\, N .
$$
Then
\se\att
\begin{eqnarray}
\nonumber
\lefteqn{
\lim_{T \to \i}\,  \frac{1}{2T}
\int_{-T}^T dt \int_0^{1/c} d\mu(\sigma) | f(s + c)|^2 \ = 
} \\
& &
\lim_{T \to \i}\,  \frac{1}{2T}
\int_{-T}^T dt \int_0^{1/c} d\mu(\sigma) \Big| \sum_{n=n_0}^{N^\prime}
 a_n n^{-(s + c)} \Big|^2  \ \ + \ O(\vare) .
\label{eqb5}
\end{eqnarray}
As
$$
\lim_{T \to \i}\,  \frac{1}{2T}
\int_{-T}^T dt\  n^{-(\s + it)} m^{-(\s -it)} \= \delta_{m\,n} \, n^{-2\s}
,
$$
we get from (\ref{eqb5}) that
\se\att
\begin{eqnarray}
\nonumber
\lefteqn{
\lim_{T \to \i}\,  \frac{1}{2T}
\int_{-T}^T dt \int_0^{1/c} d\mu(\sigma) | f(s + c)|^2 \ =
} \\
& &
\sum_{n=n_0}^{N^\prime} | a_n|^2
\int_0^{ 1/c} d\mu(\sigma)  n^{-2 \s - 2c}
\ + \ O(\vare) 
\label{eqb6}
\end{eqnarray}
for all $N^\prime \geq N$. Taking
the limit in (\ref{eqb6}) as $c$ decreases to $0$, and noting that
$\vare$ can be made arbitrarily small for $N$ large enough, 
we get 
$$
\sum_{n=n_0}^\i |a_n|^2 w_n
\= \lim_{c \to 0^+} \, \lim_{T \to \i}\,  \frac{1}{2T}
\int_{-T}^T dt \int_0^{1/c} d\mu(\sigma) | f(s + c)|^2 .
$$
But this is the same limit as (\ref{eqb4}).
\ep
\vs
Note that the integrals 
$$
\frac{1}{2T}
\int_{-T}^T dt  | f(\s + it)|^2
$$
are monotonically
decreasing as a function of $\s$, so if 
$\mu(\{0\}) = 0$, the monotone convergence theorem yields:
\begin{cor}
\label{corb1}
Assume the hypotheses of Theorem~\ref{thmb1}, and also that 
$\mu(\{0\}) = 0$. Then
\be
\label{eqb18}
\sum_{n=n_0}^\i |a_n|^2 w_n
\= \lim_{T \to \i}\,  \frac{1}{2T}
\int_{-T}^T dt \int_0^{\i} d\mu(\sigma) | f(s )|^2 .
\ee
\end{cor}
\bs
By a \emph{multiplier} of $\hw$
we mean a function $\phi$ with the property that $\phi f$ is in 
$\hw$ for every $f$ in $\hw$. It follows from the closed graph
theorem that for any multiplier $\phi$, the operator of
multiplication by $\phi$, which we denote $M_\phi$, is bounded.
It is somewhat surprising that, 
although the spaces $\hw$ consist of functions analytic in $\Oh$,
the multipliers are somehow forced to extend to be analytic on all
of $\Oo$. For the case $\mu = \mu_0$, the following theorem is due
to Hedenmalm, Lindqvist and Seip \cite{hls}.
\bt
\label{thmb2}
Let $\mu$ satisfy (\ref{eq11}) and (\ref{eq12}), and let $w_n$
be given by (\ref{eqb3}). Then the 
multiplier algebra of $\hw$ is isometrically isomorphic to 
$H^\i(\Oo) \cap \D$, where the norm on $H^\i(\Oo) \cap \D$
is the supremum of the absolute value on $\Oo$.
\et
\bp
It is clear that any multiplier $\phi$ must be in $\D$, just
by considering $\phi(s) \cdot {n_0}^{-s}$.
We shall prove the theorem in two parts:

\noindent
(A) Show that if $\phi \, \in \, \D$ has $\s_b = 0$, then 
\be
\label{eqb7}
\| M_\phi \| \= \| \phi \|_{\Oo} .
\ee
\noindent
(B) Show that if $\phi$ is a multiplier of $\hw$, then $\phi$ is
analytic and bounded in $\Oo$.
\vs
Proof of (A): Suppose $\phi(s) \= \sum_{n=1}^\i b_n n^{-s}$ is
bounded in all half-planes strictly smaller than $\Oo$.
Let $f(s) = \sum_{n=n_0}^N a_n n^{-s} $ be a finite sum in $\hw$.
Then $\phi f$ has $\s_b = 0$, so by Theorem~\ref{thmb1} 
\beq
\| M_\phi \, f \|^2 &\=& 
\lim_{c \to 0^+} \, \lim_{T \to \i}\,  \frac{1}{2T}
\int_{-T}^T dt \int_0^{\i} d\mu(\sigma) | \phi(\s+c) \, f(\s + c)|^2
\\
&\leq & \| \phi \|_{\Oo}^2 \ \| f \|^2 .
\eeq
So if $\| \phi \|_{\Oo}$ is finite, then $M_\phi$ is  bounded
on the dense subspace of $\hw$ consisting of finite sums, so
extends by continuity to be a multiplier of the whole space.

We must show that $\| M_\phi \|$ equals $\| \phi \|_{\Oo}$.
So let us assume that $\| M_\phi \| =1$ and $\| \phi \|_{\Oo} > 1$,
and derive a contradiction. (We are not assuming that 
$\| \phi \|_{\Oo}$ is necessarily finite).

For each $\s > 0$, let 
$$
N_\s \= \sup_{t} | \phi( \s + it)| .
$$
By the Phragm\'en-Lindel\"of theorem, $N_\s$ is a strictly 
decreasing
function of $\sigma$. Indeed, for $\s$ very large, $N_\s$ tends to
$|b_1|$ which is less than or equal to $\| M_\phi \| < \| \phi
\|_{\Oo} $, so the conclusion follows by applying the
Phragm\'en-Lindel\"of theorem to the function $e^{\vare s} \phi(s)$ for
an appropriate choice of $\vare$ on a vertical strip.

Moreover, in each half-plane $\O_c$ for $c > 0$, the Dirichlet
series of $\phi$ converges uniformly to $\phi$
by Bohr's theorem, 
so by Theorem~\ref{thma2}, $\phi$ is uniformly almost periodic in
$\O_c$.
Therefore there exists $\vare_1, \vare_2, \vare_3, 
\vare_4 >0 $ so that, for large enough $T$,
\be
\label{eqb8}
\Big| \{ t \, : \, | \phi (\s + it) | > 1 + \vare_1,\ -T \leq t \leq
T \} \Big| \ \geq \ \vare_2 (2T) \qquad \forall \ \vare_3 \leq \s \leq 
\vare_3 + \vare_4 .
\ee
(Just find some open set where $|\phi|$ is bigger than one, and by
uniform almost periodicty, this will recur as one moves up in the
imaginary direction).
Since multiplication by $\phi$ is a contraction, so is multiplication
by $\phi^j$ for any positive integer $j$. Therefore
$$
\| \phi^j(s) \, n_0^{-s} \|^2 \ \leq \
w_{n_0} \qquad \forall \, j \, \in \, {\mathbb N} .
$$
So by Theorem~\ref{thmb1},
we conclude that
\se\att
\begin{eqnarray}
\nonumber
w_{n_0} &\ \geq \ &
\lim_{T \to \i}\,  \frac{1}{2T}
\int_{-T}^T dt \int_0^{ \vare_4} d\mu(\sigma) | \phi^j(\s + \vare_3)|^2
n_0^{-2(\s + \vare_3)} \\
&\geq& \vare_2 \, (1 + \vare_1)^{2j}\,  n_0^{-4 \vare_3}\,
\mu([0,\vare_4]) .
\label{eqb9}
\end{eqnarray}
As the right-hand side of (\ref{eqb9}) tends to infinity with $j$,
we get a contradiction.

\vs
Proof of (B):  
Let $p_j$ denote the $j^{\rm th}$ prime, and let $\N_N$ denote the
set of positive integers all of whose prime factors are in
the set $\{ p_1, \dots, p_N\}$:
$$
\N_N \= \{  p_1^{\nu_1} \cdots p_N^{\nu_N} \ : \ \nu_1, \dots ,
\nu_N \, \in \, \N \}.
$$
For every positive integer $N$, let $Q_N$
denote orthogonal projection from $\hw$ onto the closed linear
span of the functions $$
\{ n^{-s} \, : \,  n \, \in  \, \N_N, \ n \geq n_0 \} .
$$
Suppose $\phi$ is a multiplier of $\hw$.
Then we have
\be
\label{eqb10}
Q_N  M_\phi Q_N \= M_{Q_N (\phi)} Q_N \= Q_N M_\phi .
\ee
Moreover, by a truncated version of the Euler product formula, we have that if
$f(s) = \sum a_n n^{-s}$ is any function in $\hw$, then
\beq
| Q_N (f) (s) | &\= & | \sum_{n \in \N_N} a_n n^{-s} | \\
&\leq & ( \sup |a_n|)\  \prod_{j=1}^N (1-p_j^{-\s})^{-1} .
\eeq
So if the coefficients
of $f$ are bounded, then $Q_N (f)$ is a bounded function in $\O_c$
for every $c > 0$. 
As $(Q_N f) (s + \vare) = Q_N( f(s + \vare))$, if the coefficients
of $f$ are $O(n^{\vare})$ for every $\vare > 0$, we have that
$\s_b (Q_N f) \leq \vare$ for every $\vare > 0$, \ie
$\s_b (Q_N f) \leq 0$.
By (\ref{eq13}), as the weights decay more slowly than any negative
power of $n$, it follows that for every $f$ in $\hw$,
the coefficients of $f$ are indeed
$o(n^\vare)$ for every $\vare > 0$,
and {\it a fortiori} this hypothesis holds
 for every $\phi$ in the multiplier algebra of $\hw$
(since $\phi(s) n_0^{-s}$ is in $\hw$).

Therefore we can conclude that
$$
\s_b (Q_N \, \phi) \leq 0 \qquad \forall \ N \geq 1. 
$$
Moreover, by (\ref{eqb10}), multiplication by $Q_N \phi$ on $Q_N \hw$
is a compression of $M_\phi$, so
$$
\| M_{Q_N \, \phi} \|_{Q_N \hw} \ \leq \ \| M_\phi \|_\hw \ < \ \i
.
$$
By repeating the argument in Part (A) and estimating
$$
\| (Q_N \, \phi)^j\ 2^{-\nu s} \|
$$
for $2^\nu \geq n_0$, 
we therefore conclude that
$$
\| Q_N \, \phi \|_{\Oo} \ \leq \ \|M_\phi \|_\hw \qquad \forall \,
N.
$$
By a normal families argument, some subsequence of $Q_N \, \phi$
converges uniformly on compact subsets of $\Oo$ to some $H^\i(\Oo)$
function, $\psi $ say. On compact subsets of $\O_1$, where the
Dirichlet series for $\phi$ converges absolutely, $Q_N \, \phi$
converges uniformly to $\phi$. Therefore $\phi = \psi$, and so
$\phi$ must be bounded and analytic in all of $\Oo$.
\ep

\bs
It is a theorem of  Khintchine  and Kolmogorov that if the series
$\sum |c_n|^2 $ is finite, then almost every series
$\sum \pm c_n$ converges (see \eg \cite{lev31} for a proof). 
It follows
that if $\sum a_n n^{-s}$ is in $\hw$, then for almost every choice
of signs, $ \sum \pm a_n n^{-s}$ will converge in $\Oo$ (and
in cl$(\Oo)$ for $\h_\alpha$ with
$\alpha \geq 0$). This may help explain why the
multipliers of $\hw$ extend analytically to $\Oo$.

\section{Dirichlet-like Spaces}

Throughout this section, 
let $\mu$ be a measure satisfying
conditions (\ref{eq11}) and (\ref{eq12}) of Section~\ref{secber},
and let $w_n$ be defined by (\ref{eqb3}) for $n \geq 2$. 
Define another weight
sequence  $\w^\f$ by 
$$
w_n^\f \= (\log n)^2 \ w_n .
$$
The space $\hwf$ is exactly the set of functions
whose derivatives are 
in $\hw$, and is analogous to the Dirichlet space.
We shall prove that the multipliers of $\hwf$
are contained in the multipliers of $\hw$. One can prove a similar
result for higher order derivatives, but for simplicity we stick to
the case of a single derivative.

\bt
With notation as above, the multipliers of $\hwf$ are contractively
contained in
the multipliers of $\hw$.
\label{thmc1}
\et
\bp We shall boot-strap from the following claim.

{\sc Claim:} {\it There is a constant $K < \i$ such that,
if $\phi$ is a multiplier of $\hwf$ of norm one, and both $\phi$ and
$\phip$ have $\s_b \leq 0$, then $\phi$ is a multiplier of $\hw$
of
norm at most $K$.}
\vs
Suppose the claim were proved. Let $\psi$ be any 
multiplier of $\hwf$ of norm one. 
Then for every $N$, 
$Q_N \, \psi$ satisfies the hypotheses of the claim, so
is a multiplier of $\hw$ of norm at most $K$. By
taking the weak-star limit of a subsequence of $Q_N \, \psi$, 
we can conclude that $\psi$ is a
multiplier of $\hw$ of norm at most $K$.

To show $K$ must be $1$, assume it were greater. Then there would be a
multiplier $\phi$ of $\hwf$ of norm one, which has norm greater than
$\sqrt{K}$ as a multiplier of $\hw$. By Theorem~\ref{thmb2}, 
$$
\| M_{\phi^2} \|_{\hw} \= \| M_\phi \|_\hw^2 .
$$
Then $\phi^2$ would be a multiplier of norm one of $\hwf$, and have
norm greater than $K$ as a multiplier of $\hw$, a contradiction.
\bs
\obc
We shall prove the claim with $K = \sqrt{2}$.
Suppose the claim is false. Then there is 
some finite Dirichlet
series $$
f(s) \= \sum_{n=n_0}^N a_n n^s
$$ in $\hw$ of norm $1$ such that $\| \phi f \| > K$.
Let $$
g(s) \= \sum_{n=n_0}^N a_n \frac{1}{\log n}  n^s
$$
be the primitive of $f$, which is of norm one in $\hwf$.
Let $$
B \= \{ s \, \in \, \Oo \ : \ |\phi(s)| > 1 \} .
$$
By Theorem~\ref{thmb1}, there exists $c > 0$ such that
\be
\label{eqc3}
\lim_{T \to \i}\,  \frac{1}{2T}
\int_{-T}^T dt \int_0^{\i} d\mu(\sigma) | \phi(s + c) f(s+c) |^2
\chi_B(s+c) \ > \ K^2 - 1 .
\ee
As $\phi^j$ is a multiplier of $\hwf$ of norm at most one for every positive
integer $j$, we have
$$
1\ \geq \ \| \phi^j g \|_\hwf \= \| \phi^j f + j \phi^{j-1} \phip g
\|_\hw .
$$
Therefore
\se\att
\begin{eqnarray}
\nonumber
1 &\ \geq \ &
\lim_{T \to \i}\,  \frac{1}{2T}
\int_{-T}^T dt \int_0^{\i} d\mu(\sigma) |\phi(s+c)|^{2(j-1)} \,
| \phi(s + c) f(s+c) \\
\nonumber
&&\qquad \qquad \qquad \qquad \ +\ j  \phip(s+c) g(s+c)|^2
\chi_B(s+c) \\
\nonumber &&\\
\nonumber
&\geq&
\lim_{T \to \i}\,  \frac{1}{2T}
\int_{-T}^T dt \int_0^{\i} d\mu(\sigma) 
| \phi(s + c) f(s+c) 
\\ &&\qquad \qquad \qquad \qquad \
+\ j  \phip(s+c) g(s+c)|^2
\chi_B(s+c) 
\label{eqc4}
\end{eqnarray}
By subtracting $\phi f + \phip g$ from $\phi f + j \phip g$ and
using Minkowski's inequality on (\ref{eqc4}), we get
$$
4\  \geq \ 
(j-1)^2 \ \lim_{T \to \i}\,  \frac{1}{2T}
\int_{-T}^T dt \int_0^{\i} d\mu(\sigma) 
|  \phip(s+c) g(s+c)|^2
\chi_B(s+c) 
$$
for all $j$, and so 
the limit is zero.
Therefore by Cauchy-Schwarz, (\ref{eqc4}) becomes
$$
1 \ \geq 
\lim_{T \to \i}\,  \frac{1}{2T}
\int_{-T}^T dt \int_0^{\i} d\mu(\sigma) 
| \phi(s + c) f(s+c) |^2
\chi_B(s+c) 
.$$
This contradicts (\ref{eqc3}) if $K \geq \sqrt{2}$.
\ep
\bs
Let us say a Radon measure $\nu$ supported in cl$({\Oo})$ is an 
$\alpha$-Carleson measure if there exists some constant $C$ such that
$$
\lim_{T \to \i} \frac{1}{2T} \, \int_{|\Im s| \leq T}
|f(s)|^2 d\nu(s) \ \leq \ C \| f\|^2_{\ha}
$$
for every finite Dirichlet series $f(s) = \sum_{n=2}^N a_n n^{-s}$.
Then we have
\begin{cor}
\label{corc1}
For $0 < \alpha < 2$, the function $\phi$ is a multiplier of $\ha$
if and only if 

\noindent
(i) $\phi$ is in $\D \cap H^\i(\Oo)$.

\noindent
and

\noindent
(ii) The measure $ | \phip(s) |^2 d \mu_{\alpha-2}(\s) dt$ is 
$\alpha$-Carleson.
\end{cor}
\bp
The necessity of Condition (i) follows from Theorem~\ref{thmc1}.
For Condition (ii), observe that by Cauchy's theorem, if $\phi$
is in $\D \cap H^\i(\Oo)$, then $\s_b(\phip) \leq 0$.

The function $\phi$ is a multiplier if and only if 
$$
\| \phi f^\prime + \phip f 
\|_{\h_{\alpha -2}} \ \leq \ C \|f \|_\ha
$$
for every finite Dirichlet series $f$. If $\phi$ satisfies Condition
(i), then $$
\| \phi f^\prime \|_{\h_{\alpha -2}} \ \leq \  \|\phi \|_\Oo \, \| f
\|_\ha .
$$
So such a $\phi$ is a multiplier if and only if 
$$
\| \phip f \|^2_{\h_{\alpha-2}} \=
\lim_{T \to \i} \frac{1}{2T} \int_{-T}^T \int_0^\i
| \phip(s) |^2 | f(s)|^2 d\mu_{\alpha-2}(\sigma) dt
$$
is controlled by $\| f \|^2_\ha$, \ie if and only if
Condition (ii) holds.
\ep
\bs
Of course, Corollary~\ref{corc1} does not answer the question of what
the multipliers of $\ha$ are. For this one needs a characterization of
those functions $\phi$ in $\D$ that satisfy Condition (ii), as
D.~Stegenga did for $\K_\alpha$ in \cite{ste80}.
We do not know what this condition should be.

\section{A Space of Dirichlet series with the Pick property}
\label{secpick}

Let $\h$ be a Hilbert function space on a set $X$ with reproducing
kernel $k$. We say $\h$ has the \emph{Pick property} if,
given any distinct points $\l_1, \dots , \l_N$ in $X$ and any
complex numbers $z_1, \dots, z_n$, then a necessary and sufficient
condition for the existence of a function $\phi$ in 
the closed unit ball of the multiplier algebra of $\h$ 
that has the value $z_i$ at each $\l_i$ is that
the matrix
$$
\Big[ k(\l_i, \l_j)\ ( 1 - z_i \bar z_j)  \Big]_{i,j=1}^N
$$
be positive semi-definite.
We say $\h$ has the \emph{complete Pick property} if, for any
positive integer $s$, any distinct points $\l_1, \dots , \l_N$ in
$X$ and any
$s$-by-$s$ matrices $Z_1, \dots, Z_n$,
then a necessary and sufficient
condition for the existence of a function $\phi$ in
the closed unit ball of the multiplier algebra of $\h \otimes \C^s$
that has the value $Z_i$ at each $\l_i$ is that
the $Ns$-by-$Ns$ matrix
$$
\Big[k(\l_i, \l_j)\  ( I - Z_i Z_j^\ast)  \Big]_{i,j=1}^N
$$
be positive semi-definite.
See \cite{ampi} for a treatment of complete Pick kernels.

\bs

For every integer $n \geq 2$, let $F(n)$ be the number of ways $n$
can be factored, where the order matters. Let $F(1) = 1$. Then
the following identity holds \cite[1.2.15]{tit86}
\be
\label{eqd1}
\sum_{n=1}^\i \frac{F(n)}{n^s} \= \frac{1}{2 - \zeta(s)} .
\ee
For the rest of this section, we shall fix
\be
\label{eqd15}
w_n \= \frac{1}{F(n)} ,
\ee
and consider the space $\hw$ (with $n_0 = 1$).
The kernel function for $\hw$ is then
\be
\label{eqd2}
k(s,u) \= \frac{1}{2 - \zeta(s+ \bar u)}
. \ee
As the reciprocal of $k$ has
only one positive square, it follows from
the M\raise.45ex\hbox{c}Cullough-Quiggin theorem
(see \cite{mccul92}, \cite{qui93}, \cite{agmc_cnp}),
that $k$ is a complete Pick kernel.

Let $\rho = 0.86\dots$ be the unique real number in $(1/2, \i)$
satisfying
$\zeta( 2 \rho) = 2$.
Then the space $\hw$ is analytic in the half-plane $\O_\rho$, and
in no larger domain because the kernel function blows up on the
boundary of $\O_\rho$.

As $k$ is a complete Pick kernel, there is a realization formula
for the multipliers of $\hw$.
\bt
\label{thmd05}
Let $w_n$ be defined by (\ref{eqd15}). Then $\phi$ is in the unit
ball of the multiplier algebra of $\hw$ if and only if there is an
auxiliary Hilbert space $\LL$, 
and a unitary
$U : \C \oplus \LL  \to \C \oplus \ell^2\otimes \LL $ such that,
writing
$U$ as
$$
\label{eqhb1}
U \=
\bordermatrix{&\C &\LL  \cr
\C &A & B \cr
\ell^2\otimes \LL &C  & D} ,
$$
we have 
\be
\label{eqd17}
\phi(s) \= A + B (I - E_s D)^{-1} E_s C .
\ee
Here $\ell^2 $ is the space of square summable sequences
$(c_n)_{n=2}^\i$, and $E_s : \ell^2\otimes \LL \, \to \, \LL$
is the contractive linear operator defined for $s \in \O_\rho$
on elementary tensors by
$$
E_s \ :\ (c_n)_{n=2}^\i \, \otimes \xi \ \mapsto \ \left(\sum_{n=2}^\i
c_n n^{-s} \right) \, \xi .
$$
\et
Theorem~\ref{thmd05} is a special case of the realization
formula for general complete Pick kernels (see \cite{btv},
\cite{amti00}, \cite{tom00}, \cite{ampi}), and we omit the proof.

\bs

Notice that $\phi(s)$ defined by (\ref{eqd17})
is guaranteed to be holomorphic only if $\| E_s \| = \sqrt{
\zeta(2\s) -1}$ is less than $1$, \ie if $s$ is in $\O_\rho$.
It turns out that
$\O_\rho$ is indeed the maximal domain of analyticity of the multiplier
algebra (which we shall denote $\m(\hw)$).

We shall prove this assertion via interpolating sequences.
Let $\{ \l_i \}_{i=1}^\i$ be a sequence in $\O_\rho$. We say 
$\{ \l_i \}_{i=1}^\i$ is an \emph{interpolating sequence for $\hw$}
if the linear operator
$$
\Lambda \ : \ f \ \mapsto \ \left( \frac{f(\l_i)}{\| k_{\l_i} \|} \, 
\right)_{i=1}^\i
$$
maps $\hw$ to $l^2$ surjectively. We say
$\{ \l_i \}_{i=1}^\i$ is an \emph{interpolating sequence for
$\m(\hw)$} if, for any bounded sequence $(w_i)_{i=1}^\i$, there is
a function $\phi$ in $\m(\hw)$ with $\phi(\l_i) = w_i$.

The sequence 
$\{ \l_i \}_{i=1}^\i$ is an interpolating sequence for $\hw$
if and only if the associated normalized Gram matrix
$$
G_{ij} \= \frac{\la k_{\l_j}, k_{\l_i} \ra}{\| k_{\l_j}\| \, \| k_{\l_i}
\|} 
$$
is bounded above and below (for a proof of this assertion, see \eg
\cite{nik} or \cite{ampi}). 
Therefore, every sequence 
$\{ \l_i \}_{i=1}^\i$ that tends to the line $\{ \s = \rho \}$ has
a subsequence that is an interpolating sequence for $\hw$: just
choose the $n^{\rm th}$ element in the subsequence so that the
first $n-1$ entries in the $n^{\rm th}$
row and $n^{\rm th}$ column of $G$ are smaller in modulus than
$2^{-n}$, and then the off-diagonal entries will have
Hilbert-Schmidt norm smaller than $1$.

It is a theorem of D.~Marshall and C.~Sundberg  \cite{marsun}
that, for any space with the Pick property, the interpolating
sequences for the Hilbert space and for its multiplier algebra
coincide.  Coupling this with the previous observation, we conclude
that for every boundary point of $\O_\rho$, there is a
sequence that converges to that boundary point and
that is interpolating for $\m(\hw)$. Choosing $w_i$ to alternate
$1$ and $-1$, we conclude that there is a function in $\m(\hw)$
that cannot be continued analytically to any neighborhood of that
boundary point.
Indeed we can do more. Provided $\Re(\l_i)$ converges to $\rho$
quickly enough, the sequence $\{\l_i\}$ will be interpolating for
$\hw$ regardless of how $\Im(\l_i)$ is chosen. So we can choose a
sequence that accumulates on the whole line $\{ \Re(s) \, = \,
\rho\}$, and conclude that there is a single function 
$\phi$ in $\m(\hw)$ that has $\O_\rho$ as its maximal domain of
analyticity.
Thus we have proved:
\bt
With notation as above,
the largest common domain of analyticity of the multipliers of $\hw$
is $\O_\rho$.
\label{thmd1}
\et
\bs
Is there a better way to describe the interpolating sequences for
$\m(\hw)$ (or $\hw$) than requiring that the normalized Grammian
be bounded above and below?
\begin{conjecture}
Let $w_n$ be given by (\ref{eqd15}).
Then  the sequence
$\{ \l_i \}_{i=1}^\i$ is interpolating for $\m(\hw)$
if and only if

\noindent
(i) There exists some constant $C$ such that, for any $i \neq j$,
there is a function $\phi$ in $\m(\hw)$ of norm less than $C$
with $\phi(\l_i) = 0$ and $\phi(\l_j) = 1$.

\noindent
and

\noindent
(ii) There exists some constant $C$ such that, for any $f$ in
$\hw$,
$$
\sum_{i=1}^\i \frac{|f(\l_i)|^2}{\| k_{\l_i} \|^2} \
\leq \ C  \, \| f \|_{\hw}^2  .
$$
\end{conjecture}

Condition (i) is clearly necessay, and condition (ii) is equivalent
to requiring that the normalized Grammian be bounded. So the
question is whether weak separation (condition (i)) and the
Carleson measure condition (ii) suffice. The conjecture is a
special case of Question 9.57 in \cite{ampi}.

\bibliography{references}

\end{document}